\input amstex
\documentstyle{amsppt}

\magnification=\magstephalf
\pagewidth{16true cm}
\pageheight{24.5true cm}
\frenchspacing
\parskip4pt

\font\Sans=cmss10

\font\tenfr=eufm10 
\newfam\frfam\textfont\frfam=\tenfr 
\font\eightfr=eufm8\scriptfont\frfam=\eightfr
\font\fivefr=eufm5\scriptscriptfont\frfam=\fivefr
\def\frak{\fam\frfam\tenfr}

\def\Aut{\hbox{\Sans Aut \hskip .5pt}}
\def\Hol{\hbox{\Sans Hol \hskip .5pt}}
\def\Id{\hbox{\Sans Id \hskip .5pt}}
\def\aut{\hbox{{\frak aut} \hskip .5pt}}

\topmatter

\title 
A note on a theorem of H. Cartan 
\endtitle

\author 
Jos\'e M. Isidro 
\endauthor

\affil 
Universidad de Santiago de Compostela 
\endaffil

\address 
Facultad de Matem\'aticas, 
Universidad de Santiago,
15706 Santiago de Compostela, SPAIN 
\endaddress

\email jmisidro\@zmat.usc.es \endemail


\thanks 
Supported by Comisi\'on Hispano-H\'ungara de Cooperaci\'on Cient\'{\i}fica
y Tecnol\'ogica 
\endthanks

\keywords 
Holomorphic automorphisms, Lie groups.
\endkeywords

\subjclass 
48 G 20,  72 H51 
\endsubjclass
 
\abstract 
We prove that if $D\subset \Bbb C ^n$ 
is a bounded domain with real analytic boundary, and $D$ is either 
pseudoconvex or it satisfies condition R, then the the compact
open  topology in $\Aut (D)$, the group of holomorphic
automorphisms of $D$ is the  the topology of uniform convergence
on $D$.
\endabstract

\endtopmatter

\document 

\head{0} 
Introduction 
\endhead
Let $D$ be a bounded domain in $\Bbb C ^n$ with $n\geq 1$, and let 
$\Aut (D)$ denote the group of holomorphic automorphisms of $D$. 
It was shown by H. Cartan \cite{2, Chap. 9 th. 4} that $\Aut (D)$, endowed with the 
compact open topology (that is, the topology of local uniform 
convergence on $D$), is a real Lie group whose Lie algebra $\aut (D)$ 
consists of all complete holomorphic vector fields $X\colon D\to \Bbb C ^n$,  
and that the natural action 
$$\Aut (D) \times D\to D, \qquad (f, z)\mapsto f(z)$$
is real analytic in joint variables. This result is no longer valid 
if $\Bbb C ^n$ is replaced by an arbitrary complex Banach space $E$; however, 
it is still valid if we assume that $D\subset E$ is a bounded symmetric 
domain (\cite{4}, \cite{5}). 

It is remarkable that no assumptions on the smoothness of the boundary 
$\partial D$ is needed either in the $\Bbb C ^n$ setting or in the 
infinite dimensional case. This leads naturally to the question of 
whether every $f$ in $\Aut (D)$ extends holomorphically beyond $\partial D$. 
Affirmative answers to this question have been found, in the $\Bbb C ^n$ 
setting, by S. Chen S. Jang in \cite{1, th.1.2} assuming that $D$ has 
real analytic boundary and that either $D$ is pseudoconvex or that condition 
(R) holds on $D$. By condition (R) we mean that the Bergmann projector  
(that is, the orthogonal projection from $L^2(D)$ onto the the closed subspace 
${\Cal O}^2(D)$ of square integrable holomorphic functions) maps 
${\Cal C}^{\infty}(\overline D)$ continuously into itself). In the 
infinite dimensional context, an affirmative answer has been established in 
\cite{6} for bounded circular domains with no restrictions on the boundary. 

In turn, the possibility of extending every element $f$ in $\Aut (D)$ to 
a neighbourhood of $\overline D$ gives new information about the topology of 
$\Aut (D)$. Indeed, it has been proved in \cite{6} that for bounded 
circular domains $D$ in a Banach space $E$, the topology on $\Aut (D)$ of 
local uniform convergence over $D$ is actually the topology of uniform 
convergence on $D$. In this note we use the results in \cite{1} to prove 
the following theorem:

{\bf Theorem} Let $D\subset \Bbb C ^n$, where $n\geq 1$, be a bounded 
domain with real analytic boundary, and suppose that either $D$ 
is pseudoconvex or it satisfies condition (R). Then the 
compact open topology on $\Aut (D)$ is the same as the topology of uniform 
convergence over $D$. 

For a domain $D\subset \Bbb C ^n$, we let $\Hol (D, \Bbb C ^n)$ denote the vector space 
of all holomorphic mappings $h\colon D\to \Bbb C ^n$, and the subset $\Hol (D)$ 
consists of 
the mappings such that $f(D)\subset D$. A set $S\subset D$  is a vanishing set for $D$ if
the relations
$f\in \Hol (D,\Bbb C ^n)$ and $f_{\vert  S}=0$ entail $f=0$. We say that $\partial D$ is an
algebraically determining set  for $\Aut (D)$ if the relations $f, g\in \Aut (D)$ and $f=g$
on $\partial D$ entail 
$f=g$. Similarly, the expression $\partial D$ is a topologically determining set 
for $\Aut (D)$ means that whenever $(f_n)_{n\in \Bbb N}$ and $f$ are respectively a sequence
and  an element in $\Aut (D)$ such that $f_n \to f$ uniformly over $\partial D$ we have 
$f_n\to f$ in $\Aut (D)$. Obvious changes give us the meaning of the expressions  
$S\subset D$ is an algebraically or topologically determining set for $\aut (D)$, see  
\cite{3}

\head{1} 
The main result.
\endhead

\proclaim{1.1 Theorem} 
Let $D\subset \Bbb C ^n$, where $n\geq 1$, be a bounded 
domain with real analytic boundary, and suppose that either $D$ 
is pseudoconvex or it satisfies condition (R). Then for every 
sequence $(f_n)_{n\in \Bbb N}$ in $\Aut (D)$ and every $f$ in $\Hol (D)$ 
the following conditions are equivalent:
\item{\rm (1)} $f\in \Aut (D)$ and there is a neighbourhood $\Omega$ 
of $\overline D$ in $\Bbb C ^n$ such that $(f_n)_{n\in \Bbb N}$ converges to 
$f$ uniformly on $\Omega$. 
\item{\rm (2)} $f\in \Aut (D)$ and $(f_n)_{n\in \Bbb N}$ converges to 
$f$ in the group $\Aut (D)$.
\item{\rm (3)} There are a non void open subset $U\subset D$ and a 
point $a\in U$ such that $(f_n(z))_{n\in \Bbb N}$ converges to $f(z)$ 
pointwise on $U$ and $f(a)\notin \partial D$. 
\endproclaim
\demo{Proof} We only need to prove $3\Longrightarrow 1$. By the identity 
principle, $U$ is a vanishing set for $D$. Also the family 
$(f_n)_{n\in \Bbb N}$ is uniformly bounded in $D$. Therefore by Vitali's 
theorem \cite{2, Chap. 1  prop. 7} the sequence $(f_n)_{n\in \Bbb N}$ converges 
to $f$ uniformly on each 
compact subset $K\subset D$. Since $f(a)\notin \partial D$, we have 
$f\in \Aut (D)$ by Cartan's theorem \cite{2, Chap. 5  th. 4}. On the other hand,  
since $D$ either is pseudoconvex or it satisfies condition (R), by \cite{1, th. 1.2}
there are a neighbourhood ${\Cal V}_f$ of $f$ in $\Aut (D)$ and a 
neighbourhood $\Omega_f$ of $\overline D$ in $\Bbb C ^ n$ such that every 
$g\in {\Cal V}_f$ extends to some $\tilde g \in \Hol (\Omega_f)$ and the 
action 
$$
{\Cal V}_f\times \Omega _f \to \Bbb C ^n , \qquad 
(g,z)\mapsto \tilde g (z)
$$
is real analytic on ${\Cal V}_f\times \Omega _f$. By restricting ourselves 
to smaller neighbourhoods ${\Cal W}\subset {\Cal V}_f$ and 
$\Omega \subset \Omega _1\subset \Omega_f$ we may assume that 
$\overline {\Cal W}$ and 
$\overline \Omega _1$ are compact and that the latter is connected. 
Since $f_n\to f$ in $\Aut (D)$, we have $f_n\in {\Cal W}$ for $n$ large 
enough, say $n\geq n_0$ and so 
$$ \sup _ { n\geq n_0 , z\in \Omega } \Vert f_n(z)\Vert \leq M<\infty ,$$
and a new application of Vitali's  theorem yields that $f_n\to f$ uniformly 
on each compact subset of $\Omega_1 $, in particular on $\Omega$. \qed
\enddemo
\proclaim{1.2 Corollary} Let the domain $D\subset\Bbb C ^n$ satisfy the 
assumptions in {\rm (1.1)}. 
Then 
\item{\rm (1)} The compact open topology on $\Aut (D)$ coincides with the 
topology of uniform convergence over $\overline D$ and with the topology of 
uniform convergence over $\partial D$.
\item{\rm (2)} $\partial D$ is a determining set for the topological group $\Aut (D)$. 
\item{\rm (3)} If $D$ is a bounded convex domain with real analytic boundary, 
then every element in $\Aut (D)$ has a fixed point in $\overline D$. 
\endproclaim
\demo{Proof} (1): It is clear from the above proof that on $\Aut (D)$ uniform 
convergence 
over the compact subsets of $D$ is the same as uniform convergence over $\overline D$,
which in turn is the same as uniform convergence over $\partial D$ as a consequence of 
the maximum modulus theorem. (2) and (3): Extend to $\overline D$ the elements in 
$\Aut (D)$ and apply respectively the maximun modulus theorem and Brower's fixed 
point theorem. \qed
\enddemo
We now consider the similar problems at the Lie algebra level. 
\proclaim{1.3 Theorem} 
Let $D\subset \Bbb C ^n$, where $n\geq 1$, 
be a bounded 
domain with real analytic boundary, and suppose that either $D$ 
is pseudoconvex or it satisfies condition (R). Then there are a neighbourhood 
${\Cal M}_0$ of $0$ in $\aut (D)$ and a neighbourhood $\Omega _0$ of 
$\overline D$ in $\Bbb C ^n$ such that every vector field $X\in {\Cal M}_0$ 
extends to some $\tilde X \in \Hol (\Omega_0, \Bbb C ^n)$ 
that is the infinitesimal generator of a local one-parameter group 
of holomorphic transformations on $\Omega_0$.  
\endproclaim
\demo{Proof} Choose a neighbourhood 
$V_{\Id}$ of 
$\Id$ in the group $\Aut (D)$ and a neighbourhood $\Omega _0$ 
of $\overline D$ in 
$\Bbb C ^n$ such that every element $f\in V_{\Id}$ extends to a holomorphic 
mapping $\tilde f$ in $\Omega _0$ and the action 
$$
{\Cal V}_{\Id}\times \Omega _0 \to \Bbb C ^n , \qquad 
(f,z)\mapsto \tilde f (z)$$
is real analytic on ${\Cal V}_{\Id}\times \Omega _0$. The exponential 
mapping $\exp$  defines a homeomorphism of a neighbourhood ${\Cal W}_{\Id}$ 
of $\Id$ in $\Aut (D)$ onto a neighbourhood of $0$ in $\aut (D)$. 
Let ${\Cal U}_0 \colon = {\Cal V}_{\Id}\cap {\Cal W}_{\Id}$ and set 
${\Cal M}_0\colon = \exp ^{-1} {\Cal U}_0$. We may assume that 
${\Cal M}_0$ is balanced hence connected. Let $X\in {\Cal M}_0$ and 
let $f_t \colon = 
\exp \,tX\in \Aut (D)$ denote the one-parameter group corresponding to 
$X$. For  small values of $t$, say $\vert t\vert <\tau$, we have 
$f_t\in {\Cal U}_0$, hence $f_t$ extends to a holomorphic mapping $f_t \in 
\Hol (\Omega _0 , \Bbb C ^n)$. By the identity principle we have 
$f_ t\circ f_s= 
f_{t+s}$ whenever $\vert t+ s\vert <\tau $. Thus $t\mapsto f_t$ is a local one-
parameter group of holomorphic transformations on $\Omega _0$ and so 
its infinitesimal generator $\tilde X= {d\over dt}_{\vert t=0}f_t$ is a 
holomorphic mapping on $\Omega _0$ that extends $X$. \qed
\enddemo 
\proclaim{1.4 Corollary} Let the domain $D\subset \Bbb C ^n$ satisfy the assumptions 
in {\rm (1.3)}. 
Then $\partial D$ is a determining set for $\aut (D)$, and for 
every  sequence $(X_n)_{n\in \Bbb N}$ in $\aut (D)$ and every function $X\in \Hol (D, \Bbb C
^n)$  the following conditions are equivalent:
\item{\rm (1)} $X\in \aut (D)$ and there is a neighbourhood $\Omega$ of 
$\overline D$ in $\Bbb C^n$ such that $X_n \to X$ uniformly on $\Omega$. 
\item{\rm (2)} $X\in \aut (D)$ and $(X_n)_{n\in \Bbb N}$ converges to 
$X$ in the Lie algebra $\aut (D)$.
\item{\rm (3)} 
$(X_n(z))_{n\in \Bbb N}$ is uniformly bounded in $D$ and $X_n\to X$ 
pointwise on some open non void subset $U\subset D$
\endproclaim
\demo{Proof} Clearly $1\Longrightarrow 2$. Suppose that (2) holds; then a repetition of 
the arguments made in the proof of (1.3) 
show that $X$ and all $X_n$ 
extend holomorphically to a suitable neighbourhood $\Omega$ of $\overline D$ and 
$X_n\to X$ uniformly over $\Omega$, hence in particular $(X_n)_{n\in \Bbb N}$ is 
uniformly bounded on $D$. Now suppose (3) holds; by Vitali's theorem we have 
$X_n\to X$ uniformly over the compact subset of $D$, hence $X\in \aut (D)$, 
and by (1.3)
all $X_n$ and $X$ extend holomorphically to some neighbourhood 
$\Omega$ of $\overline D$. A new application of Vitali's theorem gives the 
result. The other claim is now obvious. \qed
\enddemo

\enddocument